\documentclass[12pt,preprint]{elsarticle}

\usepackage{amsmath,amssymb}
\usepackage{subfigure}

\journal{Commun Nonlinear Sci Numer Simulat}

\begin{document}

\begin{frontmatter}

\title{Energy stable interior penalty discontinuous Galerkin finite element method for Cahn-Hilliard equation}

\author[bk]{B\"{u}lent Karas\"{o}zen}
\ead{bulent@metu.edu.tr}

\author[as]{Ay\c{s}e Sar{\i}ayd{\i}n Filibelio\u{g}lu \corref{cor}}
\ead{saayse@metu.edu.tr}

\author[bk]{Murat Uzunca}
\ead{uzunca@gmail.com}

\cortext[cor]{Corresponding author. Tel.: +90 312 2105610, Fax: +90 312 2102985}

\address[bk]{Department of Mathematics \& Institute of Applied Mathematics, Middle East Technical University, 06800 Ankara, Turkey}
\address[as]{Institute of Applied Mathematics, Middle East Technical University, 06800 Ankara, Turkey}

\begin{abstract}
An energy stable conservative method is developed for the Cahn--Hilliard (CH) equation with the degenerate mobility. The CH equation is discretized in space with the mass conserving symmetric interior penalty discontinuous Galerkin (SIPG) method. The  resulting  semi-discrete nonlinear system of ordinary differential equations are solved in time by the unconditionally energy stable  average vector field (AVF) method. We prove that the AVF method preserves the energy decreasing property of the CH equation. Numerical results confirm the theoretical convergence rates and the performance of the proposed approach.
\end{abstract}


\end{frontmatter}


\section{Introduction}

In this work, we consider Cahn-Hilliard (CH) equation in mixed form as a system of second order differential equations in a bounded domain $\Omega \subset \mathbb{R}^d \: (d\le 3)$ \cite{Chen13tss,Guo14esd}
\begin{subequations}\label{ch}
\begin{align}
u_t&= \nabla \cdot [\mu(u) \nabla w] , &\hbox{in} \:\: \Omega \times (0,T]  \\
w&=-\epsilon^2 \Delta u +f(u), &\hbox{in} \:\: \Omega \times (0,T] \\
u(x,0) & =   u_0,  &\hbox{in} \:\: \Omega
\end{align}
\end{subequations}
with the homogenous Neumann boundary conditions
\begin{equation*}
\frac{\partial u}{\partial n} =\mu(u)\frac{\partial w}{\partial n}=0, \:\: \hbox{on} \:\: \partial \Omega \times (0,T],
\end{equation*}
or with periodic boundary conditions, where $\Delta$ denotes the Laplace operator and $\mu(u)$ is the non-negative mobility function.
The parameter $\epsilon$ provides a measure of the width of inter-facial layer,  capturing the dominating effect of the reaction kinetics and represents the effective diffusivity.

The CH equation is the most known model for phase separation and coarsening phenomena in a melted alloy \cite{cahn58fen}. It was also used as diffuse  interface model for problems from fluid dynamics, material science and biology \cite{wu14ssc}. The variable $u$ denotes the concentration  of one of the species of the alloy and it is known as the phase state between materials.

The CH equation describes a gradient flow with energy dissipation
\begin{equation}\label{endis}
u_t = -\mu(u) \frac{\delta {\cal E} (u)}{\delta u}
\end{equation}
with respect to  the $H^{-1}(\Omega)$ inner product
$$
(u,v)_{H^{-1}} := (u,\Delta^{-1}v).
$$
In \eqref{endis}, the term $ \delta \mathcal{E}(u)/\delta u$ stands for the variational derivative of the Ginzburg-Landau energy functional ${\cal E}$ given by
\begin{equation} \label{energyfunction}
\mathcal{E}(u)=\int_{\Omega} \left( \frac{\epsilon^2}{2}|\nabla u|^2 + F(u)  \right)d\Omega,
\end{equation}
with the potential function $F(u)$ satisfying  $F'(u)=f(u)$.
In the literature two different types of potential functions are used: the convex double-well potential
\begin{equation}\label{func1}
F(u)=\frac{(1-u^2)^2}{4},
\end{equation}
and the non-convex logarithmic potential \cite{Guo14esd,wells06dgm}
\begin{equation}\label{logfunc}
F(u)=\frac{\theta }{2}\left[ u\ln u+ (1-u)\ln(1-u) \right] - \frac{\theta_c}{2} u^2
\end{equation}
with $ 0 < \theta \le \theta_c$, where $\theta_c$ is the transition temperature. For temperatures $\theta$ close to $\theta_c$, the logarithmic potential functional (\ref{logfunc}) is usually approximated by the quartic double-well potential (\ref{func1}). Other types of logarithmic potential functions can be found in \cite{barrett00fea,bartels11eca,xia07ldg}.
Both potential functions (\ref{func1}) and (\ref{logfunc}) satisfy the monotonicity and the Lipschitz continuity  conditions \cite{zee11goe}
\begin{align*}
( f(u_1) - f(u_2) ) (u_1 - u_2) &\geq -C_1( u_1 - u_2 )^2 , \\
\left|  f(u_1) - f(u_2)\right| &\leq L_{f} \left| u_1 - u_2\right|, \\
\left|  f'(u_1) - f'(u_2)\right| &\leq L_{f'} \left| u_1 - u_2\right|,
\end{align*}
for $u_1, u_2 \in \mathbb{R}^d$ with the constraints $|u_{1,2}|\le 1$ and with $C_1, L_{f}, L_{f'} \geq 0$. The two important  properties of the CH equation are the conservation of mass under Neumann/periodic boundary conditions \cite{Chen13tss,Guo14esd,wu14ssc}
\begin{equation*}
\int_{\Omega}u(t)dx=\int_{\Omega} u_0 dx,
\end{equation*}
and the decrease of total energy monotonically as $t \rightarrow \infty$
\begin{equation*}
\frac{d}{dt}\mathcal{E}(u(t))=- \int_{\Omega} \mu(u) \left|\nabla w \right|^2 d\Omega.
\end{equation*}
In short, the CH equation describes a conservative gradient flow.
The non-negative mobility function $\mu(u)$ can be constant or degenerate \cite{elliot96che}, where the  diffusion process is restricted to the interface zone, i.e. $u\in [-1,1]$ and it is zero outside. Commonly adapted versions of the degenerate mobility functions are $\mu(u)=\beta u(1-u)$ \cite{Guo14esd,wells06dgm} and $\mu(u)=\beta (1-u^2)$ \cite{barrett00fea,Guo14esd,wu14ssc} for a positive $\beta$.

From numerical point of  view, it is desirable to develop  numerical schemes which preserve the energy decreasing and mass conservative properties of the CH equation. Energy stability implies that the total energy of the fully discrete CH equation dissipates in time analogously to the continuous energy (\ref{energyfunction}). The schemes that preserve the discrete versions of the continuous energy lead to approximate solutions which behave qualitatively similar to the continuous ones. Explicit methods are not suitable for time discretization of the CH equation because they are not energy stable and require very small time steps due to the stability restrictions.

In the literature, the CH equation with the constant mobility function has been discretized in space using finite differences \cite{christlieb14has}, finite elements \cite{bartels11eca} and spectral methods \cite{Chen13tss}. Energy stable time discretization methods are based either on the convex splitting of the energy functional \cite{Guo14esd} or  by adding stabilization terms to the energy functional $F(u)$ \cite{tierra14}. On the other hand, the CH equation with degenerate mobility has been discretized by continuous finite elements \cite{barrett00fea,bartels11eca}, by local discontinuous Galerkin method \cite{Guo14esd,wu14ssc,xia07ldg}, discontinuous Galerkin method with $C^0$  elements and with mixed finite elements  \cite{wells06dgm}, finite differences \cite{kim07} and by spectral methods \cite{zhu99}.

In this work, we use the mass conserving symmetric interior penalty discontinuous Galerkin (SIPG) finite element method \cite{arnold82ipf,riviere08dgm} for the space discretization. Since the discontinuous Galerkin method uses piecewise polynomials, which are discontinuous at the interfaces, the discontinuous Galerkin (DG) approximation allows to capture the sharp gradients or singularities  for small $\epsilon$ that affect the numerical solution locally. For the time integrator, it is well known that the first order backward Euler method is energy stable, i.e., the discrete energy decreases without any restriction of the step size $\Delta t$ for very stiff gradient systems with $\epsilon \rightarrow 0$  \cite{hairer10epv}. The only second order implicit  energy stable method is the average vector field (AVF) method  which preserves the energy decreasing property of the gradient systems and the systems with Lyapunov functions \cite{celledoni12ped,hairer10epv}. We prove that the energy decreasing property of the fully discrete CH equation with the degenerate mobility is preserved using SIPG in space and AVF method in time. Numerical convergence rates and numerical experiments for two dimensional problems with double well (\ref{func1}) and logarithmic (\ref{logfunc}) potential functions demonstrate the  performance of the SIPG space discretization and the AVF time integrator.

The remainder of the paper is organized as follows. In the next section, Section~\ref{semidisc}, the SIPG discretization in space for the CH equation \eqref{ch} is described and the semi-discrete system of ODEs in matrix-vector form is introduced. The fully discrete system using the AVF method as a time integrator is given in Section~\ref{fulldisc}, and in Section~\ref{enst}, the energy decreasing property of the fully discrete system is proven. Numerical results are presented in Section~\ref{numeric} to demonstrate the accuracy of the numerical approach. The paper ends in Section~\ref{conc} with some concluding remarks.

\section{Semi-discrete formulations}\label{semidisc}

In this section, we outline the semi-discretization of the CH equation \eqref{ch} using symmetric interior penalty discontinuous Galerkin (SIPG), a type of discontinuous Galerkin (DG) methods, space discretization. To be being, the classical (continuous) weak formulation of the CH equation \eqref{ch} reads as: for a.e. $t\in (0,T]$
\begin{subequations}\label{2}
\begin{align}
(u_t, v)_{\Omega} + (\mu(u)\nabla w,\nabla v)_{\Omega} &= 0, & \forall v\in H_0^1(\Omega),\\
(w, v)_{\Omega} - (f(u),v)_{\Omega}-  \epsilon^2 (\nabla u,\nabla v)_{\Omega} & = 0, & \forall v\in H_0^1(\Omega),\\
(u(0),v)_{\Omega} &= (u_0,v)_{\Omega}, & \forall v\in H_0^1(\Omega),
\end{align}
\end{subequations}
where $(\cdot , \cdot)_{\Omega}$ denotes the usual $L^2$-inner product over the domain $\Omega$. It is well-known that under certain regularity assumptions, the system \eqref{2} has a unique solution in the space $L^2(0,T;H^1(\Omega))$.

In the sequel, we briefly describe the SIPG space discretization, and then, we introduce the matrix-vector form of the semi-discrete system of the CH equation \eqref{ch} as a system of ODEs.

\subsection{Discontinuous Galerkin discretization}

In this work, we use as a space discretization method for the CH equation \eqref{ch} the SIPG method \cite{arnold82ipf,riviere08dgm} which is a type of discontinuous Galerkin finite elements method. Different from the classical finite elements methods, DG methods are suitable in the use of non-conforming grids, and requires lower regularity assumptions.

Let $\mathcal{T}_{h}$ be a family of shape regular (triangular) elements $\{K_i\}\in\mathcal{T}_h$ such that $\bar{\Omega}=\cup K_i$, $K_i\cap K_j=\emptyset$ for $K_i$, $K_j \in \mathcal{T}_{h}$,  $i\neq j$. The diameter of an element $K$ and the length of an edge $E$ are denoted by $h_K$ and $h_E$, respectively. We set the space of discontinuous test and trial functions
\begin{equation}
V_h = \left\{u \in L^{2}(\Omega) : u|_{K} \in \mathbb{P}^{q}(K) \quad \forall K \in  \mathcal{T}_{h} \right\}\not\subset H_0^1(\Omega),
\end{equation}
where $\mathbb{P}^{q}(K)$ denotes the set of all polynomials on $K \in  \mathcal{T}_{h}$ of degree at most $q$. 
We split the set of all edges $E_h$ into the set $E^{0}_{h}$ of interior edges and the set $E^{\partial}_{h}$ of boundary edges so that $E_{h}=E^{\partial}_{h} \cup E^{0}_{h}$. Let the edge $E$ be a common edge for two elements $K$ and $K^{e}$. Since the functions $u\in V_h$ are discontinuous along the inter-element boundaries, there are two traces of $u$ along $E$, denoted by $u|_{E}$ from inside $K$ and $u^{e}|_{E}$ from inside $K^{e}$. Then, the jump and average of $u$ across an interior edge $E$ are defined, respectively, by
$$
[u]=u|_{E}\mathbf{n}_{K} + u^{e}|_{E}\mathbf{n}_{K^{e}}, \quad \{u\}=\frac{1}{2}(u|_{E} + u^{e}|_{E}),
$$
where $\mathbf{n}_{K}$ and $\mathbf{n}_{K^e}$ denote the outward unit normal vector to the boundary of the elements $K$ and $K^e$ on the edge $E$, respectively. Similarly, for a vector function $\nabla u$, the jump and average across an interior edge $E$ are given by
$$
[\nabla u]=\nabla u|_{E} \cdot \mathbf{n}_{K} + \nabla u^{e}|_{E} \cdot \mathbf{n}_{K^{e}}, \quad \{\nabla u\}=\frac{1}{2}(\nabla u|_{E} + \nabla u^{e}|_{E}).
$$
On a boundary edge $E \subset \partial \Omega$, we set $\{ u\}= u$ and $[u]=u\mathbf{n}$, where $\mathbf{n}$ is the outward unit normal vector to the boundary.

Then, in space SIPG discretized semi-discrete formulation of the CH equation (\ref{ch}) with a variable mobility function $\mu :=\mu (u)$ reads as: set $u_h(0), w_h(0)\in V_h$ be the projections (orthogonal $L^2$-projections) of the initial conditions $u_0, w_0$ onto $V_h$, find $u_h(t), w_h(t)\in V_h$ such that for almost every $t\in (0,T]$ and for all $\upsilon\in V_h$, we have
\begin{subequations}\label{discretecahnhilliardcoupled2}
\begin{align}
(\partial_{t}u_h, \upsilon_h)_{\Omega} + a_{h}(\mu; w_{h},\upsilon_{h})&=0,\\
(w_h, \upsilon_h)_{\Omega} - (f(u_{h}),\upsilon_{h})_{\Omega} &= a_{h}(\epsilon^2;u_{h},\upsilon_{h}),
\end{align}
\end{subequations}
where the bi-linear (in last two argument) form $a_{h}(\kappa;w,\upsilon)$ is given by
\begin{equation}
\begin{aligned} \label{dgbilinear}
a_{h}(\kappa;w,\upsilon)&=\sum_{K \in \mathcal{T}_h}\int_K \kappa\nabla w \cdot\nabla \upsilon - \sum_{E\in E^{0}_{h}}\int_E\left\{\kappa\nabla w\right\}\cdot[\upsilon]ds \\
                     & \;  - \sum_{E\in E^{0}_{h}}\int_E\left\{\kappa\nabla \upsilon \right\}\cdot[w]+\sum_{E\in E^{0}_{h}} \frac{\sigma \kappa}{h_{E}}\int_E [w]\cdot[\upsilon]ds,
\end{aligned}
\end{equation}
where the parameter $\sigma$ is called the penalty parameter and it should be sufficiently large to ensure the stability of the SIPG discretization as described in \cite{riviere08dgm} with a lower bound depending only on the polynomial degree \cite{riviere08dgm}. In numerical experiments we take $\sigma =3q(q+1)$.

\subsection{Semi-discrete system in matrix-vector form}

Let $u_h(t)$ and $w_h(t)$ be the time dependent solutions of the semi-discrete system \eqref{discretecahnhilliardcoupled2}. Then, $u_h(t)$ and $w_h(t)$ are of the form
\begin{equation}\label{5}
u_h(t)=\sum^{N}_{m=1}\sum^{n_q}_{j=1}\xi^{m}_{j}(t) \varphi^{m}_{j}\;, \quad w_h(t)=\sum^{N}_{m=1}\sum^{n_q}_{j=1}\zeta^{m}_{j}(t) \varphi^{m}_{j},
\end{equation}
where $\varphi^{m}_{j}$ are the basis functions spanning the space $V_h$, and $\xi^{m}_{j}$ and $\zeta^{m}_{j}$ are the unknown coefficients. The number $n_q$ stands for the local dimension depending on the basis polynomial order $q$ and $N$ is the number of (triangular) elements. By substituting the expansions \eqref{5} into the system \eqref{discretecahnhilliardcoupled2} and choosing the test functions $\upsilon=\varphi^{k}_{i}$, $i=1, \ldots, n_q$, $k=1, \ldots, N$, we obtain the semi-linear systems of ordinary differential equations
\begin{equation}\label{7}
\begin{aligned}
M\xi_{t}+ A_{\mu}\zeta&=0, \\
A_{\epsilon}\xi + b(\xi )-M\zeta&=0,
\end{aligned}
\end{equation}
with the ordered unknown coefficient vectors and the basis functions
\begin{align*}
\xi = \{\xi^j\}_{j=1}^{n_q\times N} &= (\xi_{1}^{1},\ldots , \xi_{n_q}^{1},\xi_{1}^{2},\ldots , \xi_{n_q}^{2},\ldots , \xi_{1}^{N}, \ldots , \xi_{n_q}^{N})^T, \\
\zeta = \{\zeta^j\}_{j=1}^{n_q\times N}&= (\zeta_{1}^{1},\ldots , \zeta_{n_q}^{1},\zeta_{1}^{2},\ldots , \zeta_{n_q}^{2},\ldots , \zeta_{1}^{N}, \ldots , \zeta_{n_q}^{N})^T, \\
\varphi = \{\varphi^j\}_{j=1}^{n_q\times N} &= (\varphi_{1}^{1},\ldots , \varphi_{n_q}^{1},\varphi_{1}^{2},\ldots , \varphi_{n_q}^{2},\ldots , \varphi_{1}^{N}, \ldots , \varphi_{n_q}^{N})^T,
\end{align*}
In \eqref{7}, $M$ denotes the mass matrix with the entries $M_{ij}=(\varphi^{j}, \varphi^{i})_{\Omega}$, $1\leq i,j \leq n_q\times N$, $A_{\mu}$ and $A_{\epsilon}$ are the stiffness matrices with the entries $(A_\mu)_{ij}=a_h(\mu(u_h);\varphi^{j},\varphi^{i})$ and $(A_{\epsilon})_{ij}=a_h(\epsilon^2;\varphi^{j},\varphi^{i})$, $1\leq i,j \leq n_q\times N$, and $b$ is the non--linear vector of unknown coefficient vector $\xi$ with the entries $b_{i}(\xi)=(f(u_h),\varphi^{i})_{\Omega}$, $1\leq i \leq n_q\times N$.

\section{Fully discrete system}\label{fulldisc}

In this section, we give the fully discrete formulations of the CH equation \eqref{ch} through the time integration of the semi-linear system of ODEs \eqref{7}. To do this, we consider the uniform partition $0=t_0<t_1<\ldots < t_J=T$ of the time interval $[0,T]$ with the uniform time step-size $\Delta t=t_{k}-t_{k-1}$, $k=1,2, \ldots, J$.

In this work, as a time integrator we use the average vector field (AVF) method which is a structure preserving time integrator for gradient systems  \cite{celledoni12ped}. For a given system of ODEs $\dot{y} = -\nabla U(y)$, the AVF method is given by
\begin{equation}\label{avfgrad}
y_n = y_{n-1} - \Delta t \int _{0}^{1} \nabla U (\tau y_n+( 1-\tau )y_{n-1})d \tau.
\end{equation}
The AVF method possesses the energy decreasing property without restriction  of the step sizes $\Delta t$ and is second order accurate in time. It represents a modification of the implicit mid-point rule and  for quadratic potentials  $U(y)$, the AVF method reduces to the mid-point rule. Higher order variants of the AVF methods for Hamiltonian and Poisson  systems with Gauss-Legendre collocation points are given in \cite{hairer10epv}. As Gauss-Legendre Runge-Kutta methods, the AVF method and higher order versions do not have damping property for very stiff systems, whereas for discontinuous Galerkin-Petrov methods and Radau II Runge-Kutta methods, the energy decreases monotonically without restriction of the step size $\Delta t$ and the Lipschitz constant for $\nabla U(y)$. However, they require the solution of coupled system of equations, which increases the computational cost for two and three dimensional CH equations, where efficient solution techniques are required \cite{estep00uks}.

Let for $t=0$, $u_{h,0}, w_{h,0}\in V_h$ be the projections (orthogonal $L^2$-projections) of the initial conditions $u_0$, $w_0$ onto $V_h$, and let $\eta_0=(\xi_0,\zeta_0)^T$ be the corresponding initial coefficient vector each of components satisfying \eqref{5}. At a specific time $t=t_n$, we denote the coefficient vector of the approximate solutions $(u_{h,n},w_{h,n})^T$ by $\eta_n=(\xi_n,\zeta_n)^T$. Then, the application of the AVF method \eqref{avfgrad} to the semi-discrete system \eqref{discretecahnhilliardcoupled2} leads to solving for $n=0,1,\ldots , J-1$ the non-linear system of equations
\begin{equation}\label{8}
\begin{bmatrix}
M & \frac{\Delta t}{2} A_{\mu}\\
\frac{1}{2} A_{\epsilon}&-\frac{1}{2} M
\end{bmatrix}
\begin{bmatrix}
\xi_{n+1} \\
\zeta_{n+1}
\end{bmatrix}+\begin{bmatrix}
 0\\
\int^{1}_{0} b(\tau \xi_{n+1} + (1-\tau)\xi_n)d\tau
\end{bmatrix}
=
\begin{bmatrix}
M\xi_n - \frac{\Delta t}{2}A_{\mu} \zeta_n\\
\frac{1}{2}M\zeta_n - \frac{1}{2}A_{\epsilon} \xi_n
\end{bmatrix},
\end{equation}
where we used the fact that the AVF method reduces the mid-point rule for linear terms. Further, the solution dependent mobility function $\mu(u_h)$ in the bilinear form $A_{\mu}$ is computed explicitly as for continuous finite elements in \cite{barrett00fea}, i.e., on each time interval $(t_n,t_{n+1}]$ the mobility function is taken as $\mu = \mu (u_{h,n})$ where $u_{h,n}$ is the known approximate solution from the previous time step. The system \eqref{8} can be written as a residual equation $R(\eta)=(R_1(\eta),R_2(\eta))^T=0$ with
\begin{equation}\label{nonliniems}
\begin{aligned}
R_1(\eta_{n+1})&=M(\xi_{n+1} - \xi_{n}) + \frac{\Delta t}{2} A_{\mu}(\zeta_{n+1}+\zeta_n) \\
R_2(\eta_{n+1})&= \frac{1}{2}A_{\epsilon}(\xi_{n+1} + \xi_{n}) - \frac{1}{2} M(\zeta_{n+1}+\zeta_n) + \int^{1}_{0} b(\tau \xi_{n+1} + (1-\tau)\xi_n)d\tau .
\end{aligned}
\end{equation}
We solve the nonlinear system of equations \eqref{nonliniems} using the Newton's method.  Starting with an initial guess $\eta_{n+1}^{(0)}=(\xi_{n+1}^{(0)}, \zeta_{n+1}^{(0)})^T$, the $k-{th}$ Newton iteration to solve the nonlinear system of equations (\ref{nonliniems}) for the unknown vector $\eta_{n+1}=(\xi_{n+1}, \zeta_{n+1})^T$ reads as
\begin{equation}\label{newties}
Js^{(k)} = - R(\eta_{n+1}^{(k)}), \qquad
\eta_{n+1}^{(k+1)} = \eta_{n+1}^{(k)} + s^{(k)} \; , \quad k=0,1,\ldots
\end{equation}
until a user defined tolerance is satisfied. In (\ref{newties}), $s=(s_1,s_2)^T$ is the increment, and $J$ stands for the Jacobian matrix
\begin{equation}\label{jacob}
J= \left(\begin{array}{cc}
M &  \frac{\Delta t}{2}A_{\mu}\\
\frac{1}{2}A_{\epsilon} + J_b &-\frac{1}{2}M
\end{array}\right) \; , \quad J_b=\tau\int^{1}_{0} b'(\tau \xi_{n+1} + (1-\tau)\xi_n)d\tau
\end{equation}
where $b'(\tau \xi_{n+1} + (1-\tau)\xi_n)$ is the Jacobian of the nonlinear form $b(\xi)$ w.r.t. $\xi$ at $\xi =\tau \xi_{n+1} + (1-\tau)\xi_n$. At each Newton iteration, the integral term in $J_b$ is approximated by the  fourth order Gaussian quadrature formula.

\section{Energy stability}\label{enst}

It is expected that fully discrete energy stable schemes should preserve the discrete energy dissipation as their continuous parts, which leads to qualitatively better approximations.  The  continuous (in time) energy  of the semi-discrete CH equation is given by \cite{feng14aip}
\begin{equation}\label{discreteenergy}
\mathcal{E}^h(u)=\frac{\epsilon^2}{2}\left\|\nabla u\right\|^{2}_{L^2(\mathcal{T}_h)} + (F(u),1)_{\Omega}.
\end{equation}
On the other hand, the discrete DG counterpart of the continuous energy \eqref{discreteenergy} at a time $t^n=n \Delta t$ reads as
\begin{equation}\label{energy}
\begin{aligned}
\mathcal{E}^h_{DG}(u^{n}) &= \frac{\epsilon^2}{2}\left\|\nabla u^{n}\right\|^2_{L^2(\tau_h)} + (F(u^n),1)_{\Omega} \\
 & \; + \sum_{E\in E^{0}_{h}} \left( -(\{\epsilon^2\partial _n u^{n}\}, [u^{n}])_E + \frac{\sigma\epsilon^2}{2h_E}([u^{n}],[u^{n}])_E \right).
\end{aligned}
\end{equation}

In this section, we show that the AVF method applied to the semi-discrete system \eqref{discretecahnhilliardcoupled2} is energy stable through the discrete energy \eqref{energy}. Applying the AVF method to the semi-discrete system \eqref{discretecahnhilliardcoupled2}, we obtain for any $\nu , \vartheta\in V_h$
\begin{equation}\label{fullyDGavf}
\begin{aligned}
(u^{n+1} - u^{n}, \nu )_{\Omega} + \frac{\Delta t}{2} a_h(\mu (u^n); w^{n+1}+w^n,\nu ) &= 0, \\
\left(\frac{w^{n+1}+w^n}{2},\vartheta\right)_{\Omega} - \int_0^1(f(\tau u^{n+1} + (1-\tau)u^n), \vartheta)_{\Omega} d\tau &= \frac{1}{2}a_h(\epsilon^2; u^{n+1}+u^n,\vartheta).
\end{aligned}
\end{equation}
Taking $\nu =(w^{n+1}+w^n)/2$ and $\vartheta=u^{n+1}-u^n$ in \eqref{fullyDGavf}, we obtain
\begin{eqnarray*}
\left(u^{n+1} - u^{n}, \frac{w^{n+1}+w^n}{2}\right)_{\Omega} + \frac{\Delta t}{4} a_h(\mu (u^n); w^{n+1}+w^n,w^{n+1}+w^n) &=& 0,\\
\left (\frac{w^{n+1}+w^n}{2},u^{n+1}-u^n\right )_{\Omega} - \int_0^1(f(\tau u^{n+1} + (1-\tau)u^n), u^{n+1}-u^n)_{\Omega} d\tau &=& \\
 \frac{1}{2}a_h(\epsilon^2; u^{n+1}+u^n,u^{n+1}-u^n).
\end{eqnarray*}
By using the identity $(a+b, a-b)_{\Omega}= (a^2 - b^2 ,1)_{\Omega}$ and the bi-linearity of $a_h$, we get
\begin{align}
\left (u^{n+1} - u^{n}, \frac{w^{n+1}+w^n}{2}\right)_{\Omega} + \frac{\Delta t}{4} a_h(\mu(u^n); w^{n+1}+w^n,w^{n+1}+w^n) &= 0, \label{eq1avf} \\
\left (\frac{w^{n+1}+w^n}{2},u^{n+1}-u^n\right)_{\Omega} - \int_0^1(f(\tau u^{n+1} + (1-\tau)u^n), u^{n+1}-u^n)_{\Omega} d\tau &= \nonumber \\
 \frac{1}{2}a_h(\epsilon^2; u^{n+1},u^{n+1}) -\frac{1}{2}a_h(\epsilon^2; u^n,u^n). \label{eq2avf}
\end{align}
Using the Taylor expansions of the terms $F (u^{n})$ and $F (u^{n+1})$, and neglecting the higher order terms, we get
\begin{align*}
F(u^n) &\approx F(\tau u^{n+1} + (1-\tau)u^n) - f(\tau u^{n+1} + (1-\tau)u^n)(\tau (u^{n+1} - u^n)),  \\
F(u^{n+1})&\approx F(\tau u^{n+1} + (1-\tau)u^n) + f(\tau u^{n+1} + (1-\tau)u^n)(1-\tau) (u^{n+1} - u^n)).
\end{align*}
Subtracting $F(u^n)$ from $F(u^{n+1})$ leads to
\begin{align}
F(u^{n+1})-F(u^{n}) &\approx f(\tau u^{n+1} + (1-\tau)u^n) (u^{n+1} - u^n) \nonumber \\
(F(u^{n+1}),1)_{\Omega}-(F(u^{n}),1)_{\Omega} &\approx ( f(\tau u^{n+1} + (1-\tau)u^n), u^{n+1} - u^n)_{\Omega} \nonumber \\
\int_0^1((F(u^{n+1}),1)_{\Omega}-(F(u^{n}),1)_{\Omega})d\tau &\approx \int_0^1( f(\tau u^{n+1} + (1-\tau)u^n), u^{n+1} - u^n)_{\Omega}d\tau \nonumber \\
(F(u^{n+1}),1)_{\Omega}-(F(u^{n}),1)_{\Omega} &\approx \int_0^1( f(\tau u^{n+1} + (1-\tau)u^n), u^{n+1} - u^n)_{\Omega}d\tau \label{tayavf}
\end{align}
We note  that the bi-linear form $a_h$ satisfies
\begin{align}
a_h(\mu(u^n); w^{n+1}+w^n,w^{n+1}+w^n) &=  \mu(u^n) \| \nabla(w^{n+1}+w^n)\|_{L^2(\Omega)}^2 \nonumber \\
                                    & \; - 2\sum_{E \in E^{0}_{h}} \int_E \mu(u^n)  \{ \nabla (w^{n+1}+w^n)\}[ w^{n+1}+w^n]ds \nonumber \\
                                    & \; +  \sum_{E \in E^{0}_{h}} \frac{\sigma \mu(u^n)}{h_E} \|[w^{n+1}+w^n]\|_{L^2(E)}^2 . \label{posa}
\end{align}
Since all the terms in \eqref{posa} are non-negative (see \cite[Sec. 2.7.1]{riviere08dgm} for positivity of edge integral term), we have $a_h(\mu(u); w^{n+1}+w^n,w^{n+1}+w^n)\geq 0$. Similarly, we have $a_h(\epsilon^2; u^{n+1},u^{n+1})\geq 0$, $a_h(\epsilon^2; u^{n},u^{n})\geq 0$. Using this identities,  subtracting \eqref{eq1avf} from \eqref{eq2avf} and substituting \eqref{tayavf}, we obtain
\begin{align*}
- \frac{\Delta t}{4} a_h(\mu(u^n); w^{n+1}+w^n,w^{n+1}+w^n) &\approx  (F(u^{n+1}), 1)_{\Omega} +\frac{1}{2}a_h(\epsilon^2; u^{n+1},u^{n+1})  \\
&  \; -  \left( (F(u^n),1)_{\Omega} + \frac{1}{2}a_h(\epsilon^2; u^{n},u^{n}) \right) \le 0,\\
\end{align*}
which implies that $\mathcal{E}^h_{DG}(u^{n+1}) \leq \mathcal{E}^h_{DG}(u^{n})$. Hence, the AVF discretized scheme is energy stable through the discrete energy \eqref{energy}.

\section{Numerical results}\label{numeric}
In this section, we present a set of numerical examples to confirm the accuracy, energy stability and mass conservation of our approach. All examples are considered in two dimensional spatial domain with constant and degenerate mobility function, double-well \eqref{func1} and logarithmic \eqref{logfunc} potential functions  under Neumann and periodic boundary conditions.

\subsection{Constant mobility and double-well potential under Neumann boundary conditions}\label{ex1}

We consider  2D CH equation with the constant mobility function $\mu(u)=1$, diffusivity $\epsilon=0.1$, and double-well potential \eqref{func1} under the homogenous Neumann boundary condition \cite{Chen13tss}. The
The exact solution in $\Omega=[-1,1]\times [-1,1]$ for $t\in [0,1]$ is given by
$$
u(x,y,t)=e^{\cos(t)}\cos(\pi x)\cos(\pi y).
$$
We add a load vector to the system so that the exact solution above solves the system. The $L^2$ errors and the numerical order of accuracy at time $t=1$ using the time step $\Delta t = 0.5 \Delta x$ for linear and quadratic DG polynomials are given in Table \ref{table1}. Both DG approximations satisfy  the $(k+1)$-th order of accuracy for $\mathbb{P}^{k}$ elements.

\begin{table}[htb]
\centering
\caption{Example \ref{ex1}: Accuracy test with constant mobility and double-well potential.\label{table1}}
\begin{tabular}{ l l  l l l }\hline
&$\Delta x$ & Dof &$L^{2}$--Error& Order\\
\hline
& 1/2 & 24 &3.347 &-\\
$\mathbb{P}^{1}$& 1/4 & 96&1.633 &1.04\\
& 1/8 & 384 &4.810E-01 &1.76\\
& 1/16 & 1536 &1.079E-01 &2.16\\
\hline
& 1/2  & 48 &6.694E-01&-\\
$\mathbb{P}^{2}$& 1/4  & 192 & 2.685E-01 &1.32\\
& 1/8 & 768 &3.376E-02 &2.99\\
& 1/16 &  3072 &3.733E-03 &3.18\\
\hline
\end{tabular}
\end{table}

\subsection{Degenerate mobility and double-well potential  under periodic boundary conditions}\label{ex2}
The next example is the CH equation  in \cite{Guo14esd} with the exact solution
$$
u(x,y,t)=e^{-2t}\sin(x)\sin(y)
$$
in the domain $\Omega=[0,2\pi] \times [ 0,2\pi]$ for $t\in [0,1]$ under periodic boundary conditions with the degenerate mobility function $\mu(u)=1-u^2$ and double-well potential \eqref{func1}. The effective diffusivity is taken as $\epsilon=1$. For first order DG polynomials we use time step $\Delta t =0.0032 \pi$ and for quadratic DG polynomials
$\Delta t =0.00032 \pi$. Order reduction is observed in Table \ref{table2} for linear and quadratic DG polynomials which is due to the nonlinearity of the degenerate mobility function.

\begin{table}[htb]
\centering
\caption{Example \ref{ex2}: Accuracy test with degenerate mobility and double-well potential.\label{table2}}
\begin{tabular}{ l l l l l}\hline
&Mesh Size  & Dof &$L^2$ error  & Order\\
\hline
& $\pi/2$  & 24 & 2.054  &-\\
$\mathbb{P}^{1}$& $\pi/4$  &96 &  5.742E-01  &1.84\\
& $\pi/8$ &384& 1.566E-01&1.87\\
& $\pi/16$ & 1536& 5.478E-02  &1.52\\
\hline
& $\pi/2$ &48 &4.342E-01&-\\
$\mathbb{P}^{2}$& $\pi/4$ &192 &1.136E-01 &1.93\\
& $\pi/8$ & 768 &1.713E-02   &2.73\\
& $\pi/16$ &  3072&4.895E-03 &1.81\\
\hline
\end{tabular}
\end{table}

\subsection{Constant mobility and double-well potential under Neumann boundary conditions:  spinodal decomposition and nucleation}\label{ex3}
We consider CH equation with constant mobility $\mu(u)=1$ and double-well potential function \eqref{func1} under homogenous Neumann boundary condition \cite{Guo14esd}. The computational domain is taken as $\Omega=[0,1]\times [0,1]$ with $\epsilon = 10^{-5}$.

This problem represents the two main separation mechanisms: spinodal decomposition and nucleation. Both mechanisms in the CH equation are defined by the initial condition $u_0(x)=\bar{u} + r$ where $\bar{u}$ is a constant and $r$ is a random number uniformly distributed on $[−0.005,0.005]$.

\begin{figure}[htb!]
\centering
\subfigure[$t=0.002$]{\includegraphics[width=0.45\textwidth] {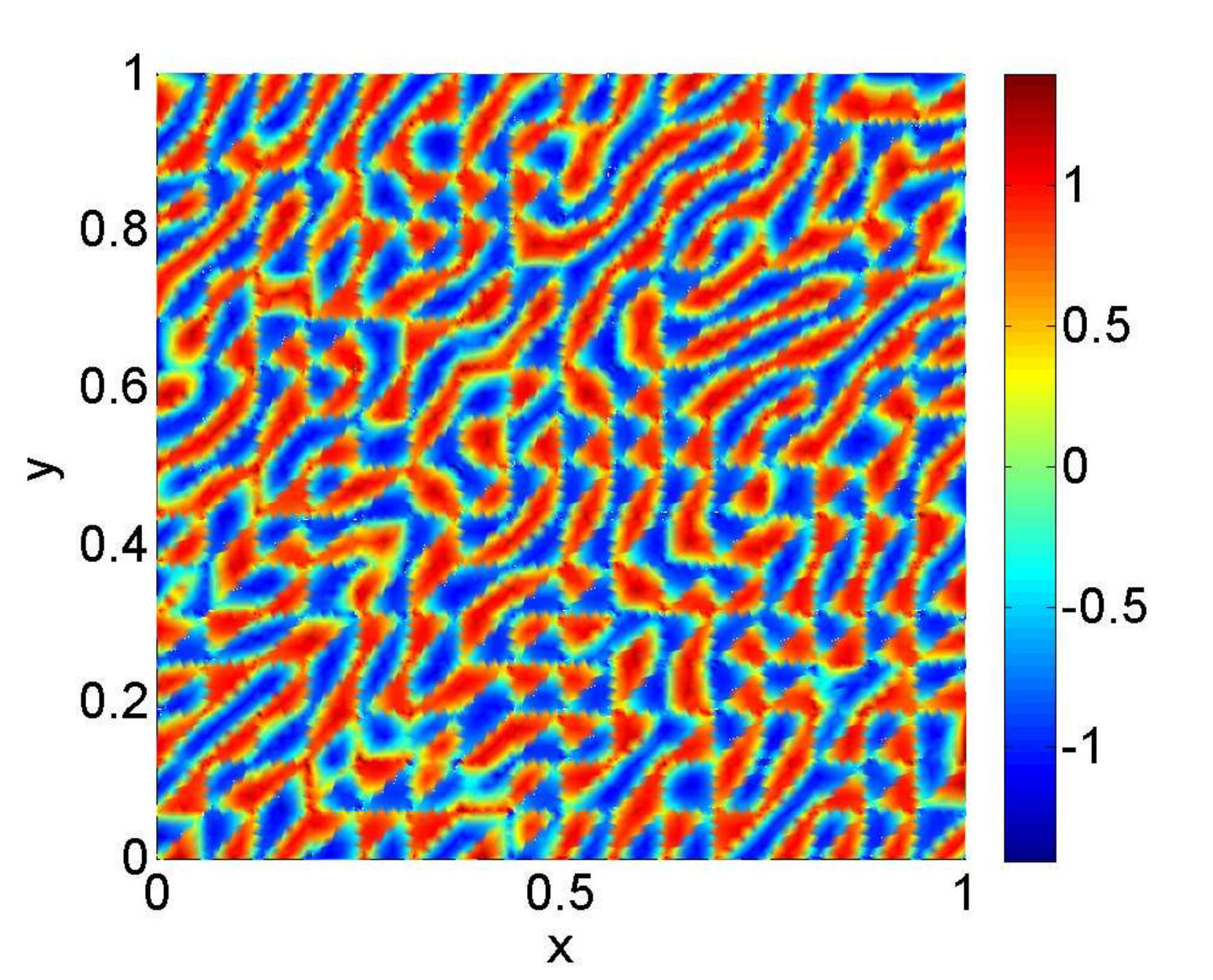}}
\subfigure[$t=0.004$]{\includegraphics[width=0.45\textwidth] {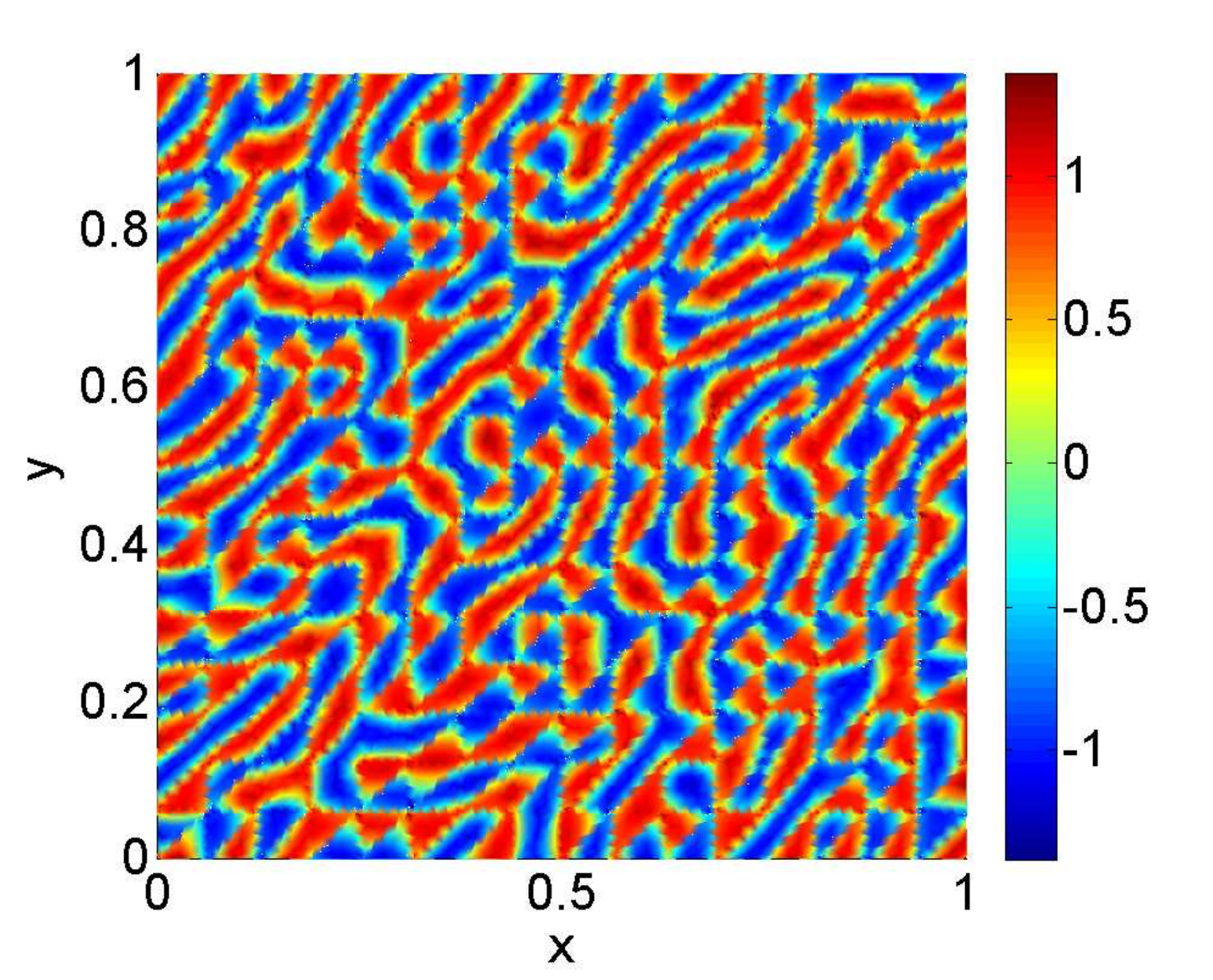}}

\subfigure[$t=0.01$]{\includegraphics[width=0.45\textwidth] {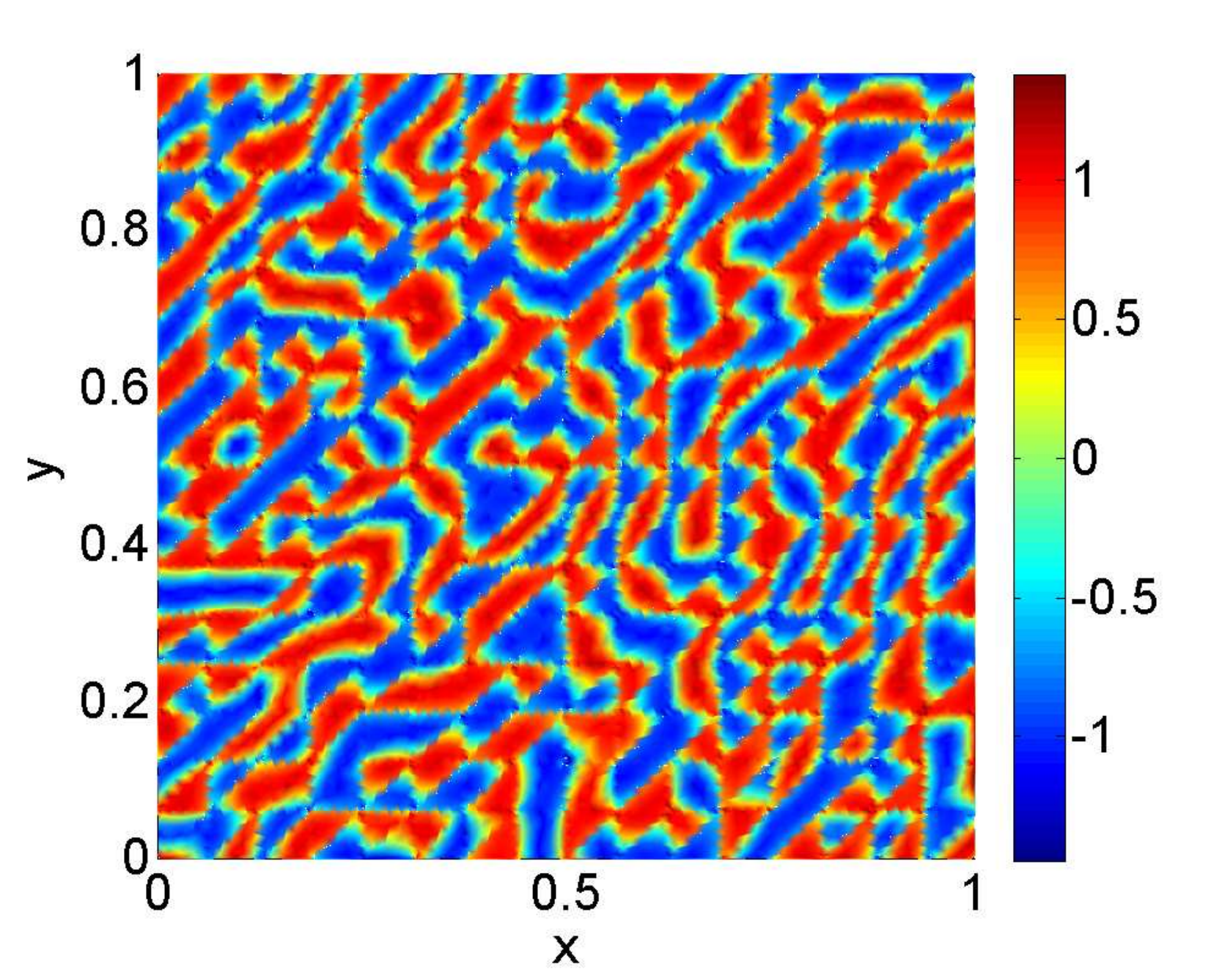}}
\subfigure[$t=0.4$]{\includegraphics[width=0.45\textwidth] {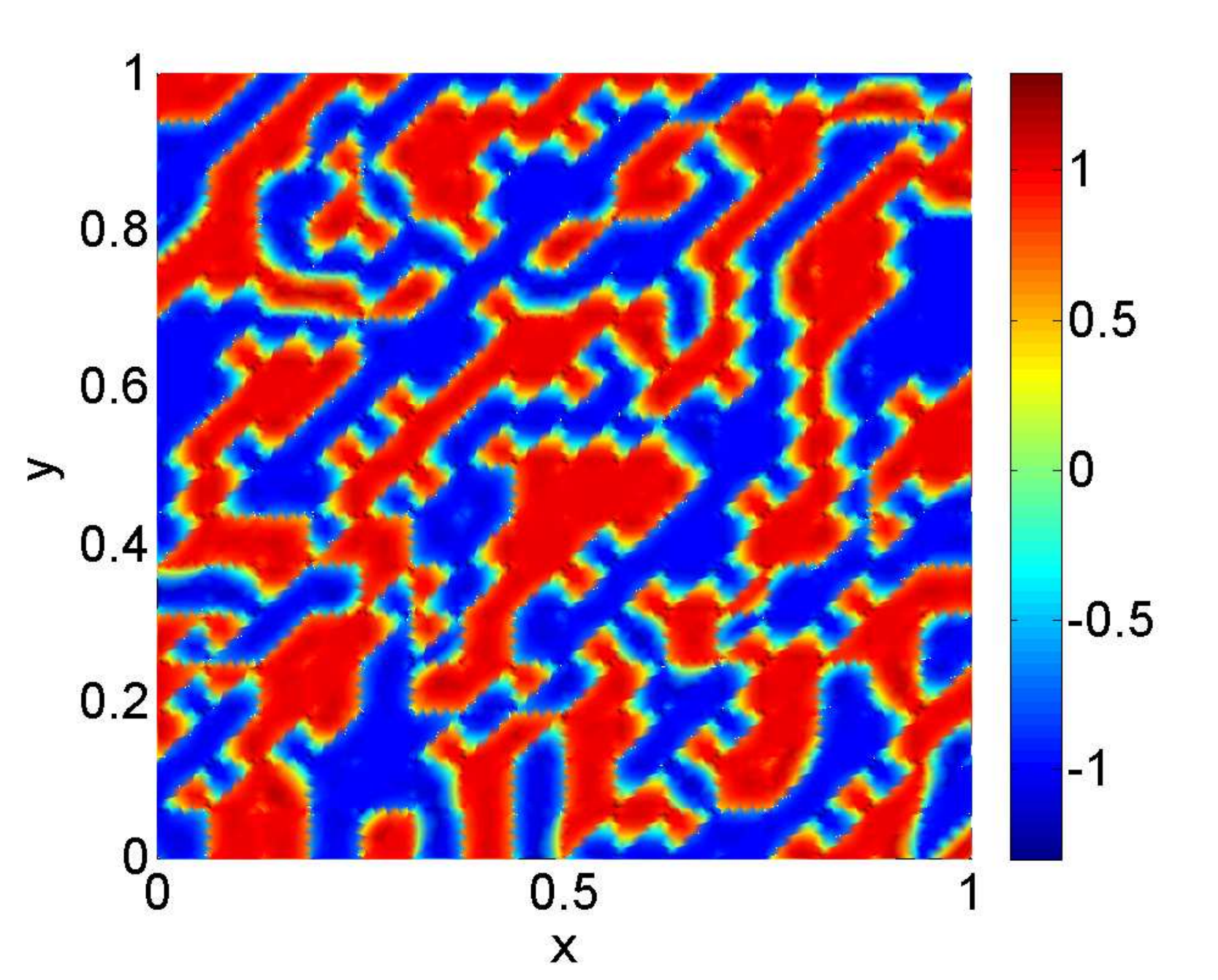}}

\subfigure[]{\includegraphics[width=0.45\textwidth] {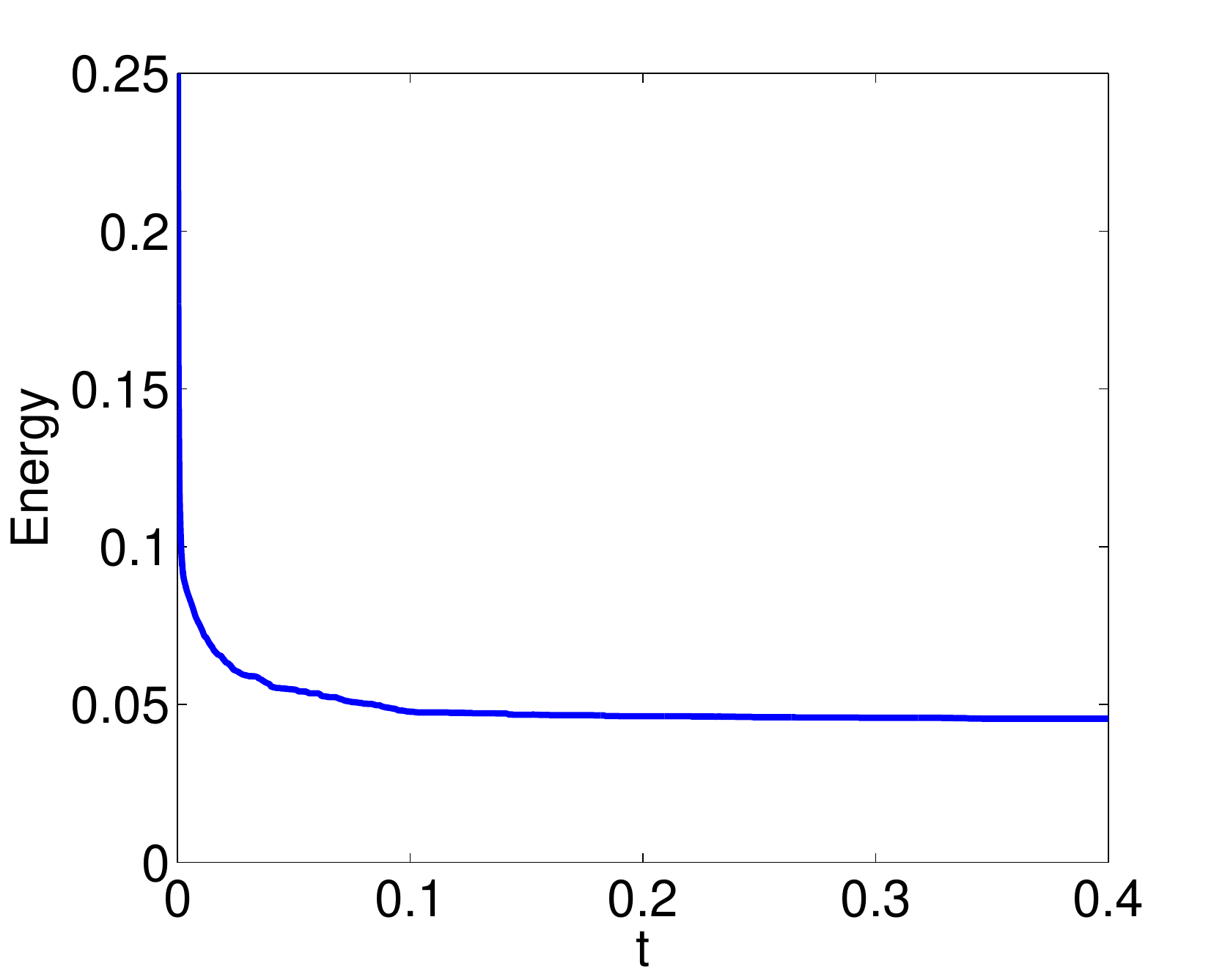}}
\subfigure[]{\includegraphics[width=0.45\textwidth] {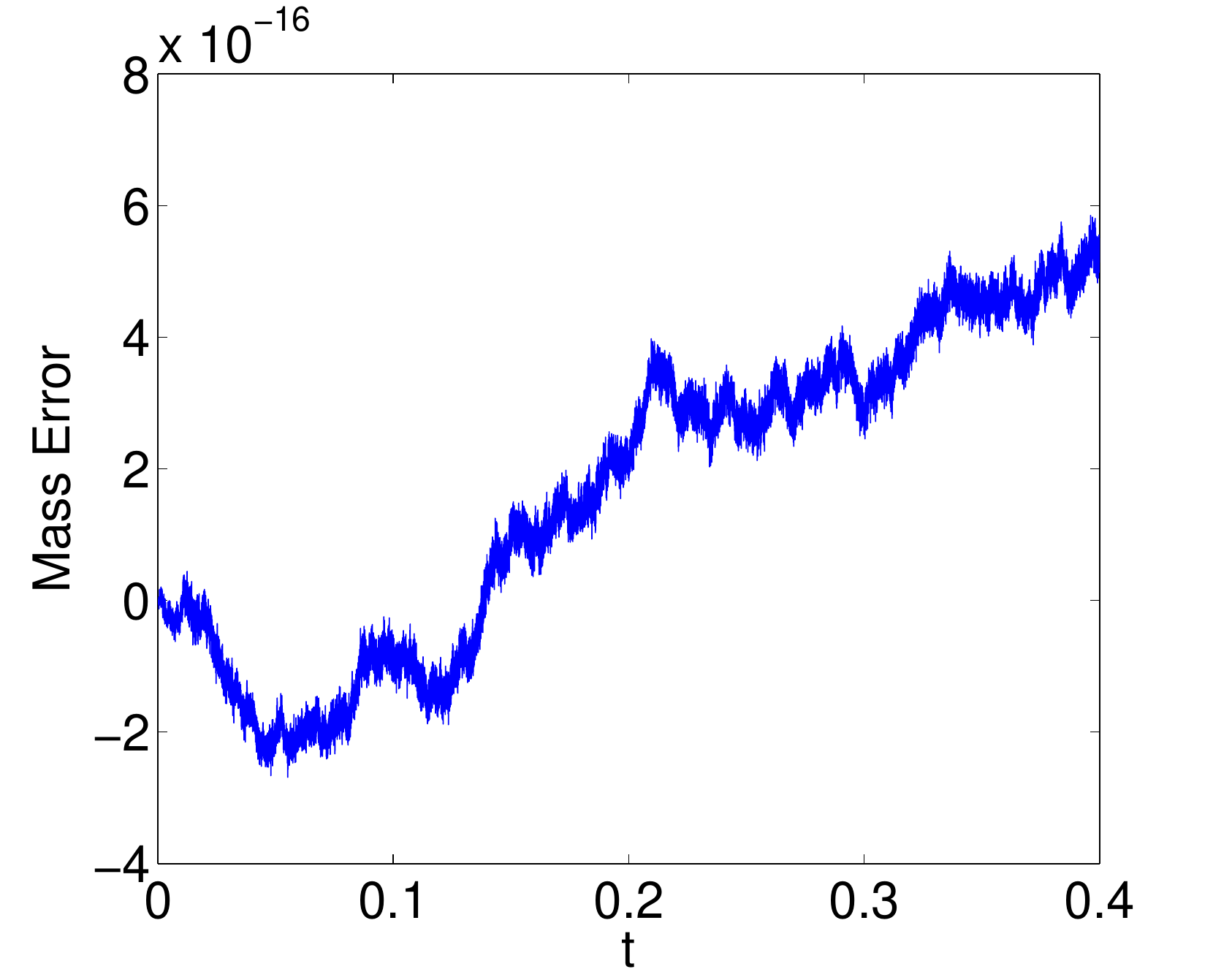}}
\caption{Example \ref{ex3}: (Top-Middle) Solutions with constant mobility and double-well potential (spinodal decomposition, $\bar{u}=0$), (Bottom) energy and mass evolutions. Linear DG elements and $\Delta t=10^{-5}$.\label{const}}
\end{figure}

\begin{figure}[htb]
\centering
\subfigure[$t=0.002$]{\includegraphics[width=0.45\textwidth] {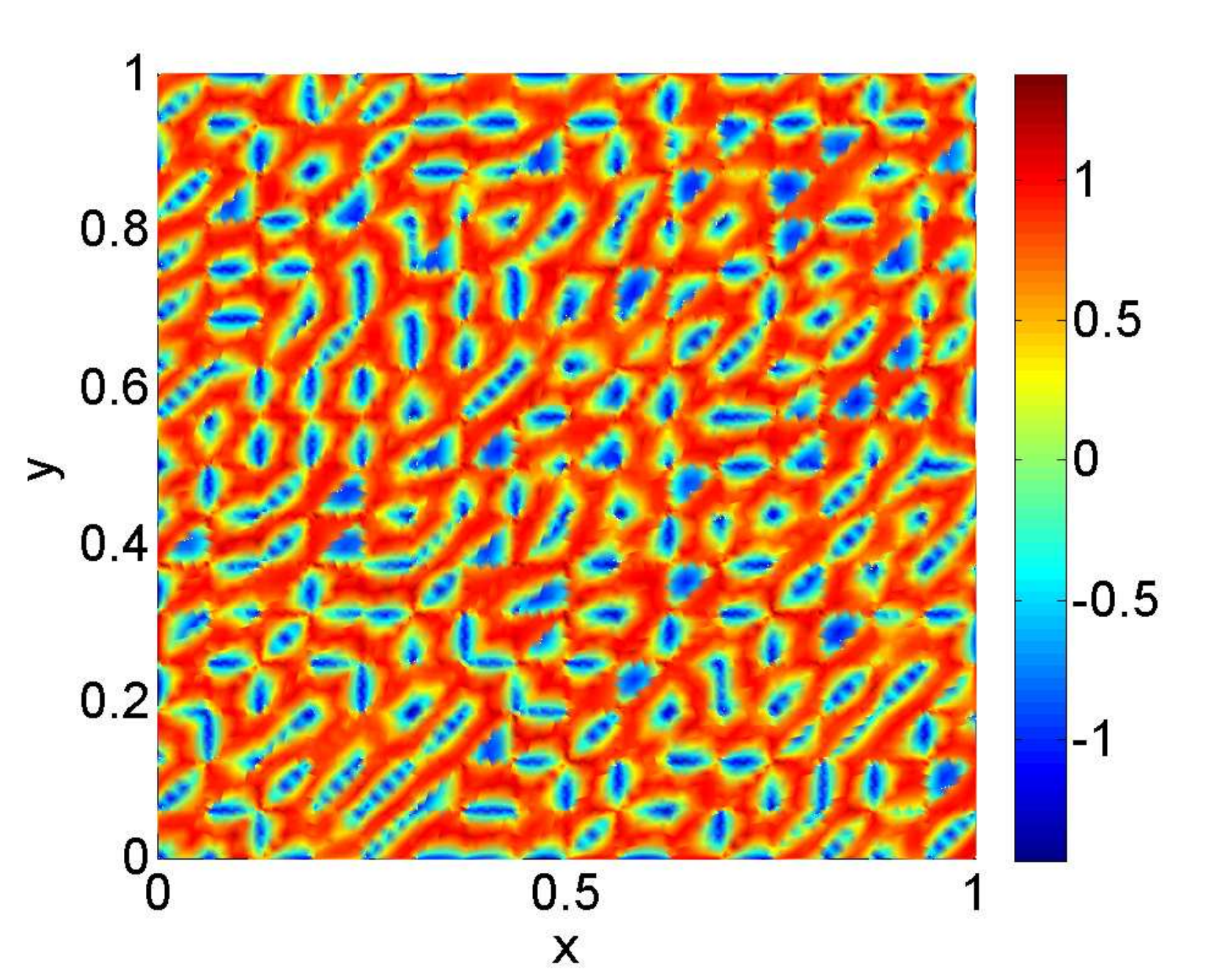}}
\subfigure[$t=0.004$]{\includegraphics[width=0.45\textwidth]{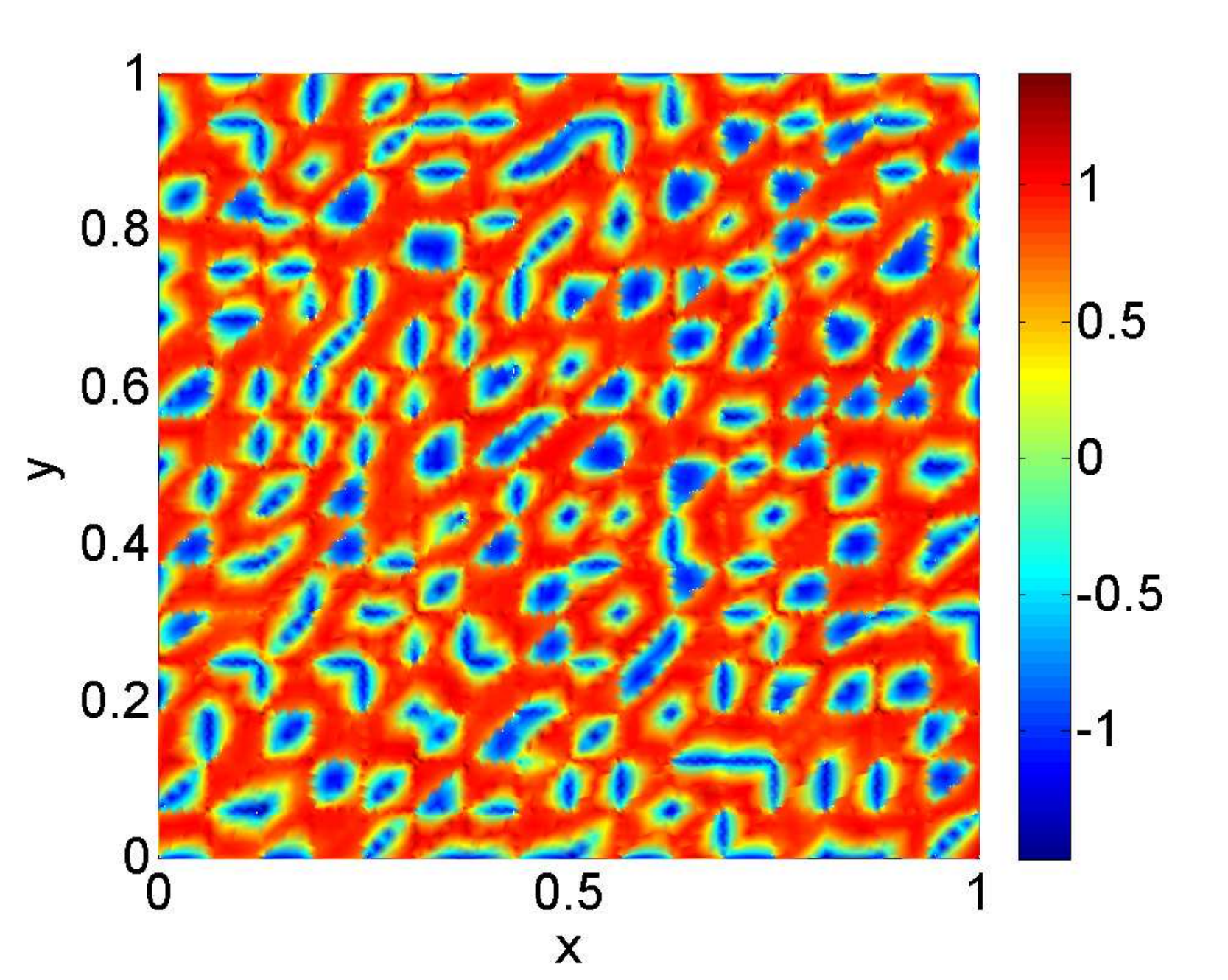}}

\subfigure[$t=0.01$]{\includegraphics[width=0.45\textwidth]{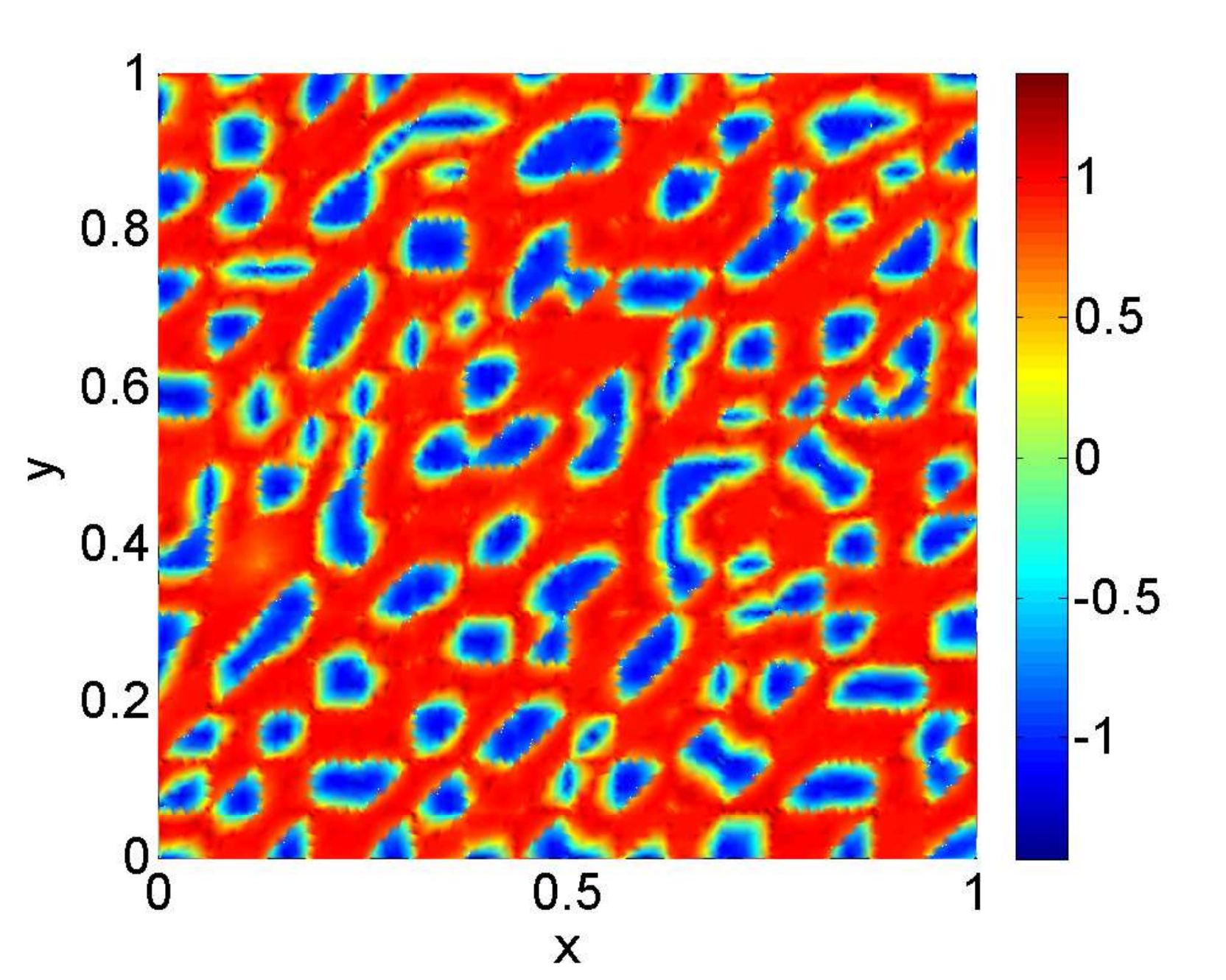}}
\subfigure[$t=0.4$]{\includegraphics[width=0.45\textwidth]{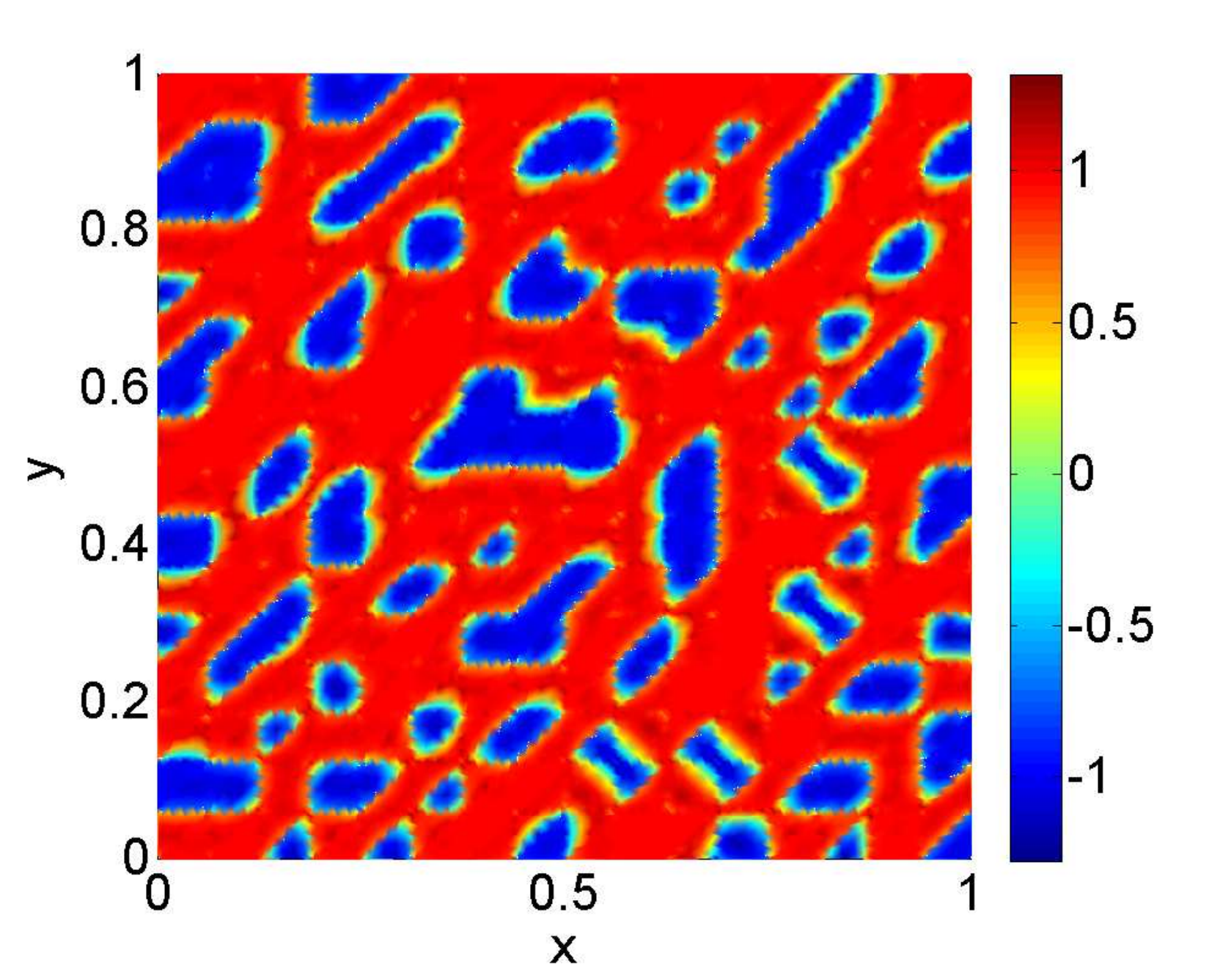}}

\subfigure[]{\includegraphics[width=0.45\textwidth] {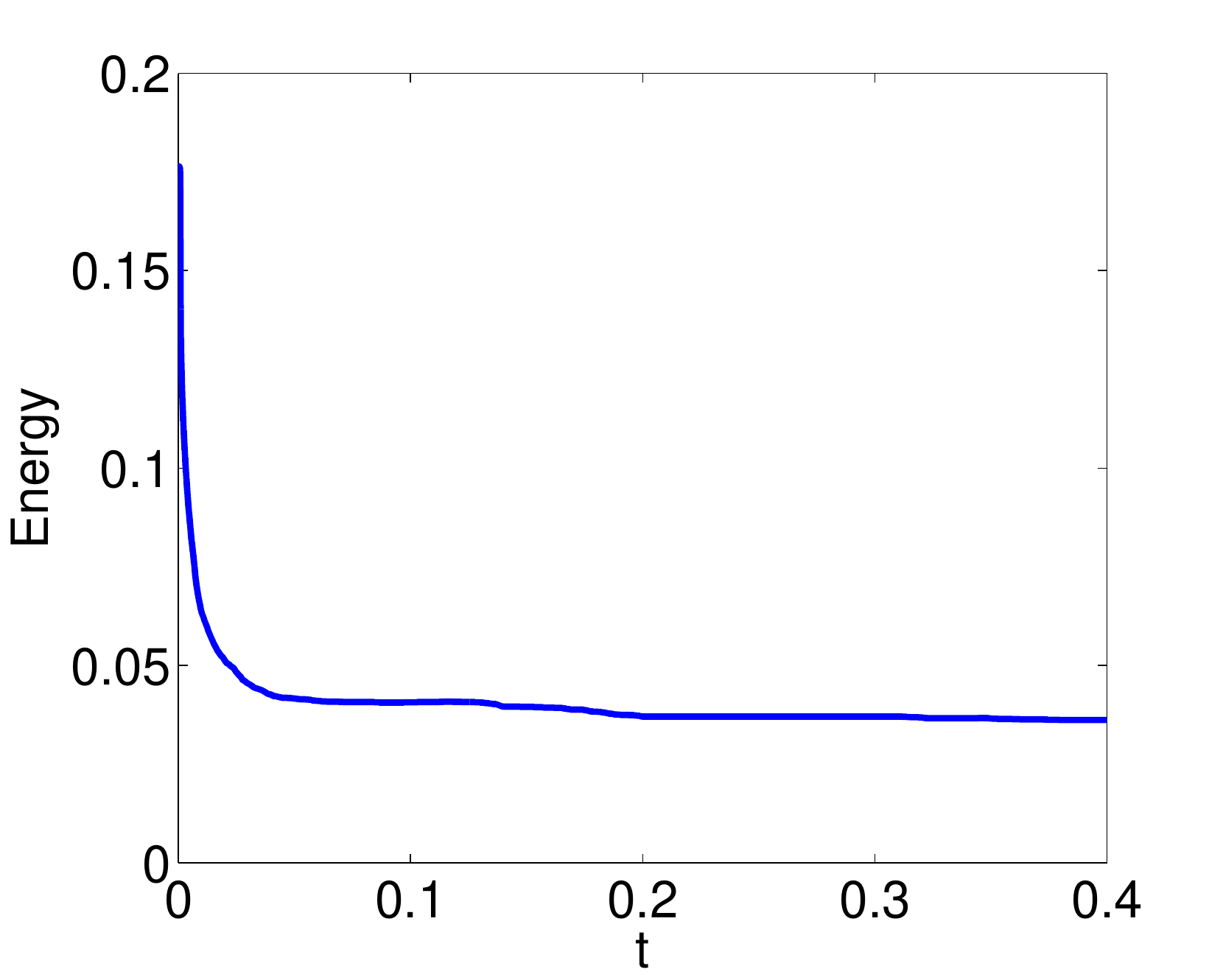}}
\subfigure[]{\includegraphics[width=0.45\textwidth] {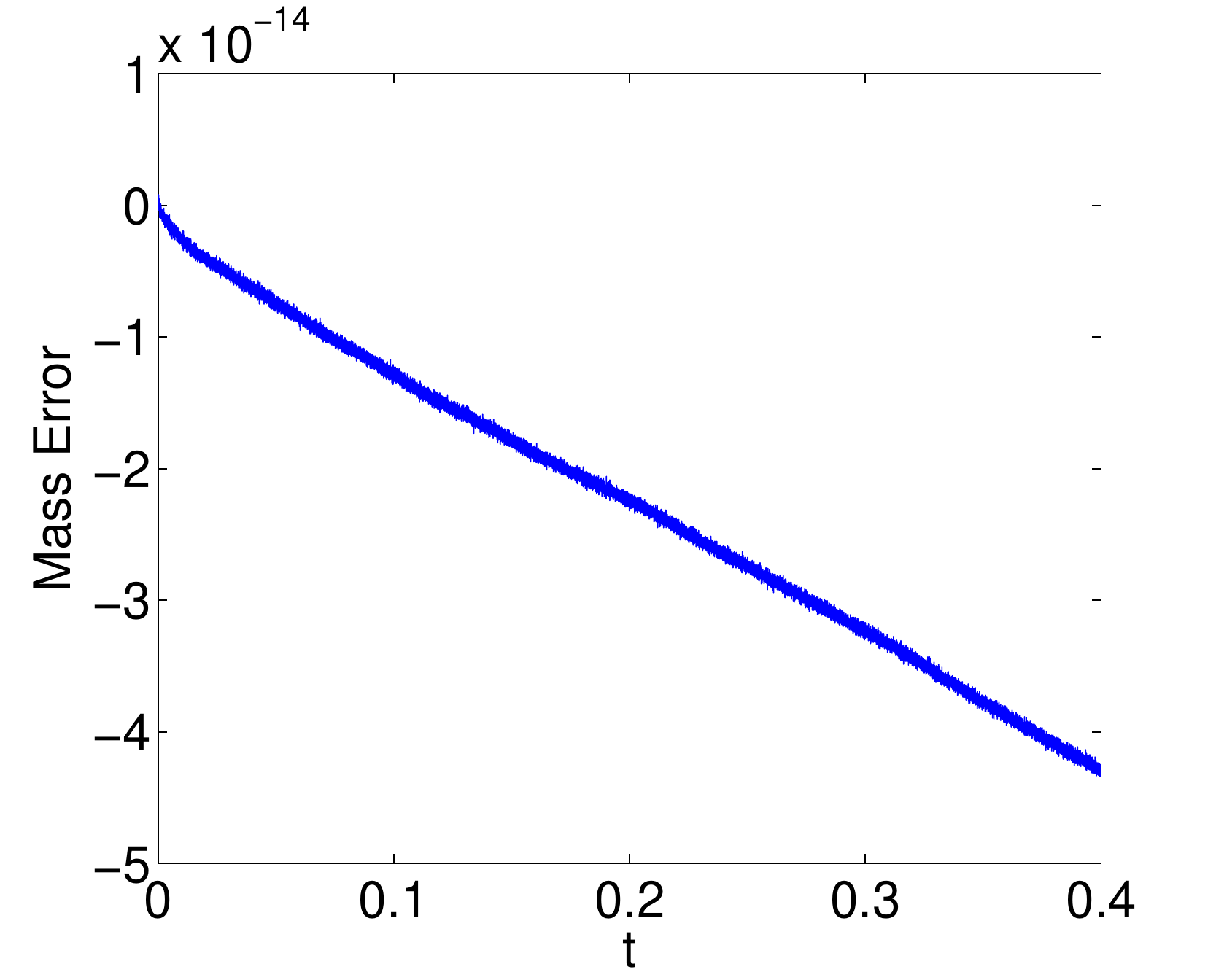}}
\caption{Example \ref{ex3}: (Top-Middle) Solutions with constant mobility and double-well potential (nucleation, $\bar{u}\neq 0$), (Bottom) energy and mass evolutions. Linear DG elements and $\Delta t=10^{-5}$. \label{const1}}
\end{figure}

For $\bar{u}=0$, spinodal decomposition will be the governing separation mechanism (see Fig.~\ref{const}). The mixture separates from a randomly perturbed homogeneous state ($\bar{u}=0$) giving rise to a striped pattern of complicated topology that coarsens over time. If we let the simulation evolve, the stationary solution would be a fully separated flow with two rectangular patches.
While $\bar{u}\neq 0$, the separation mechanism is nucleation (see Fig.~\ref{const1}). We show the solution for $\bar{u}=0.4$. In the nucleation mechanism, isolated nuclei come up from the mixture. Again, the spatial micro-structure of the mixture coarsens over time.

In both cases the discrete energy dissipates and the mass is conserved. Our results are similar to those in \cite{gomez14aeg}, where for spatial discretization local discontinuous Galerkin method and for time discretization implicit convex splitting have been used.


\subsection{Degenerate mobility and  logarithmic potential  under Neumann boundary conditions}\label{ex4}
We consider 2D CH equation with degenerate mobility function $\mu(u)=u(1-u)$ and the logarithmic potential function
$$
F(u)=3000( u\ln u + (1-u)\ln (1-u)) + 9000u(1-u)
$$
under homogenous Neumann boundary conditions  \cite{Guo14esd} with the diffusion constant $\epsilon=1$ in the domain $\Omega=[-0.5,0.5]\times[-0.5,0.5]$ for $t\in [0,0.2]$.  The initial condition is a random variation of uniform state $u=0.63$ with a change no larger than $0.05$.
\begin{figure}[htb]
\centering
\subfigure[$t=0.000002$]{\includegraphics[width=0.45\textwidth] {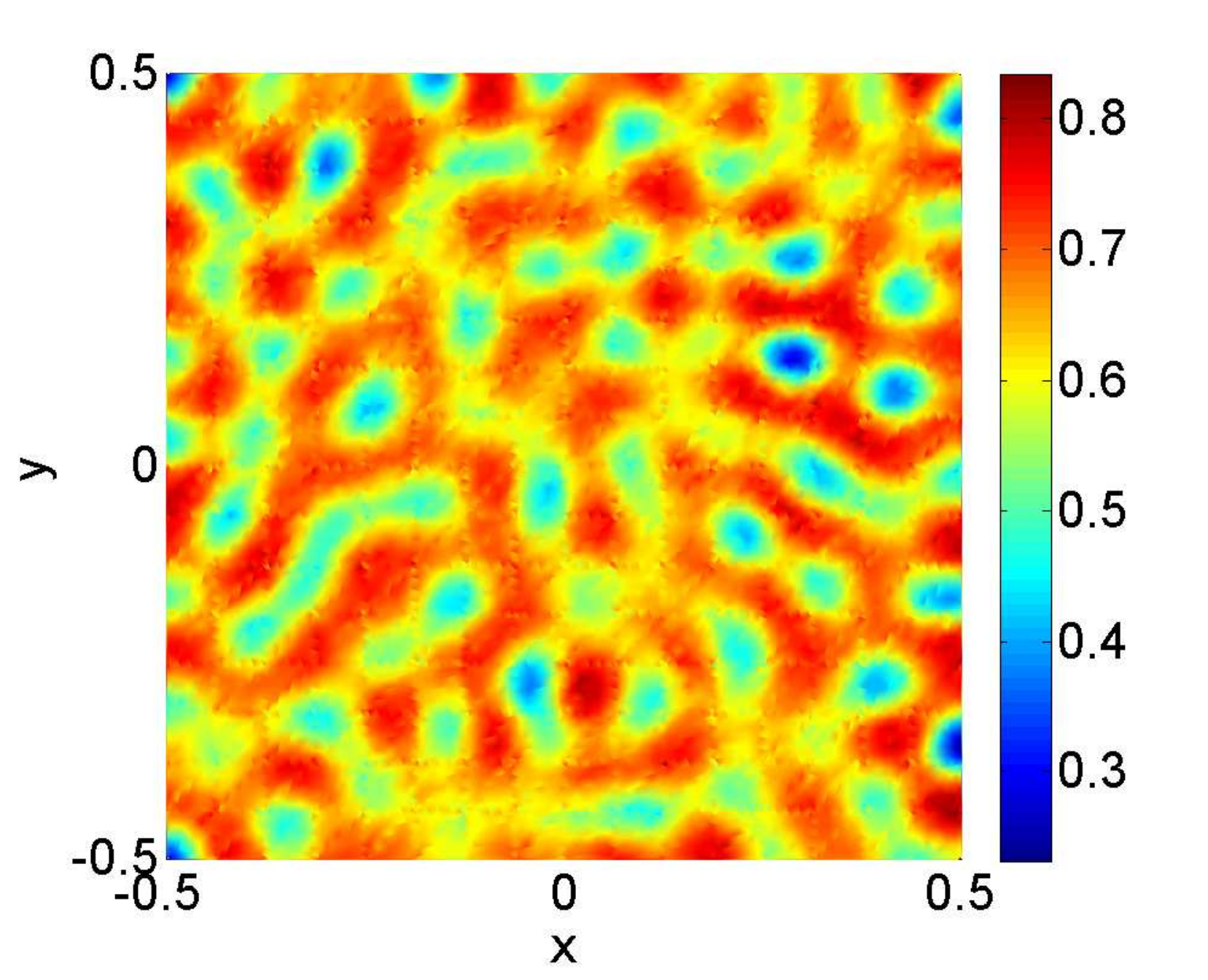}}
\subfigure[$t=0.000008$]{\includegraphics[width=0.45\textwidth] {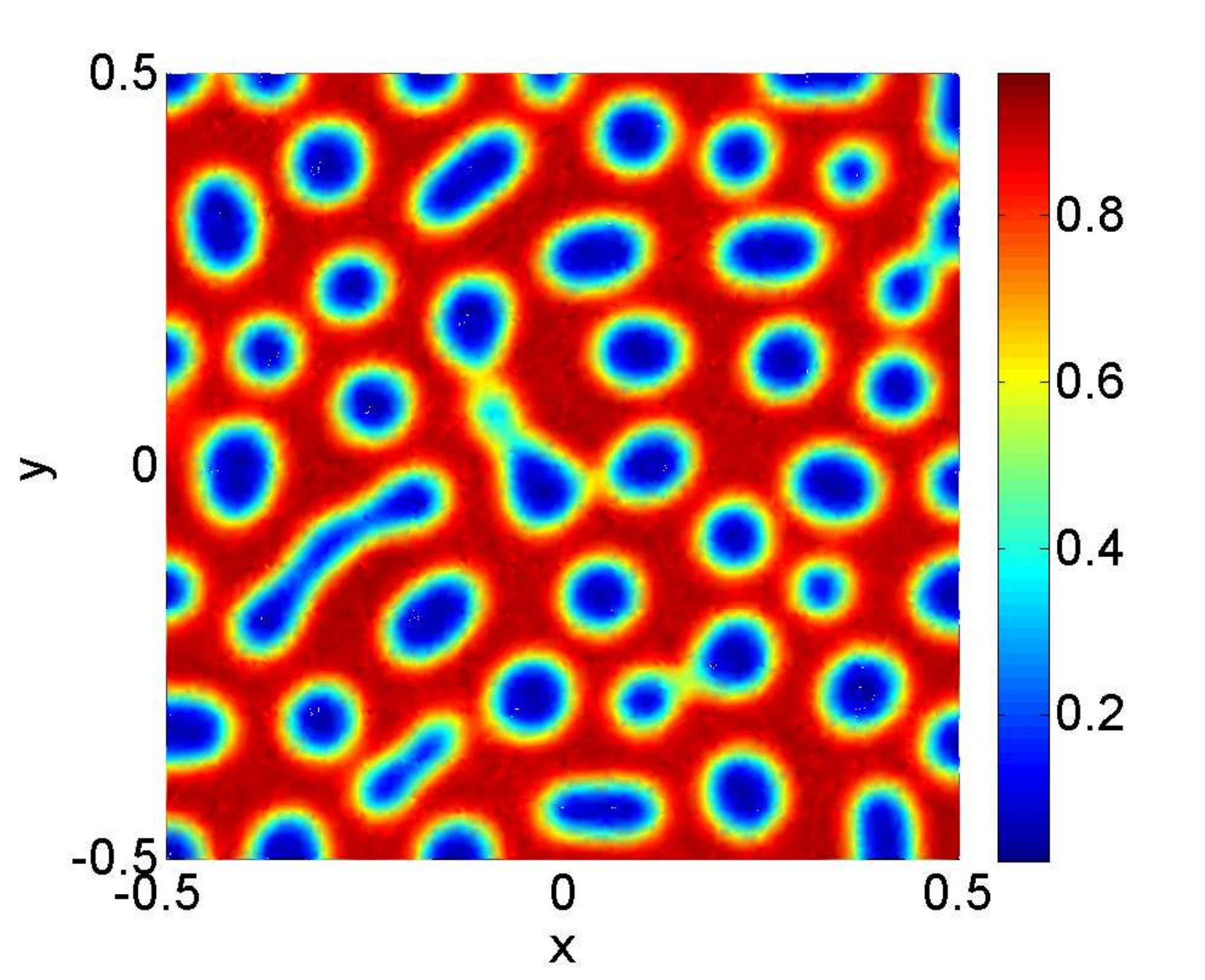}}

\subfigure[$t=0.000128$]{\includegraphics[width=0.45\textwidth] {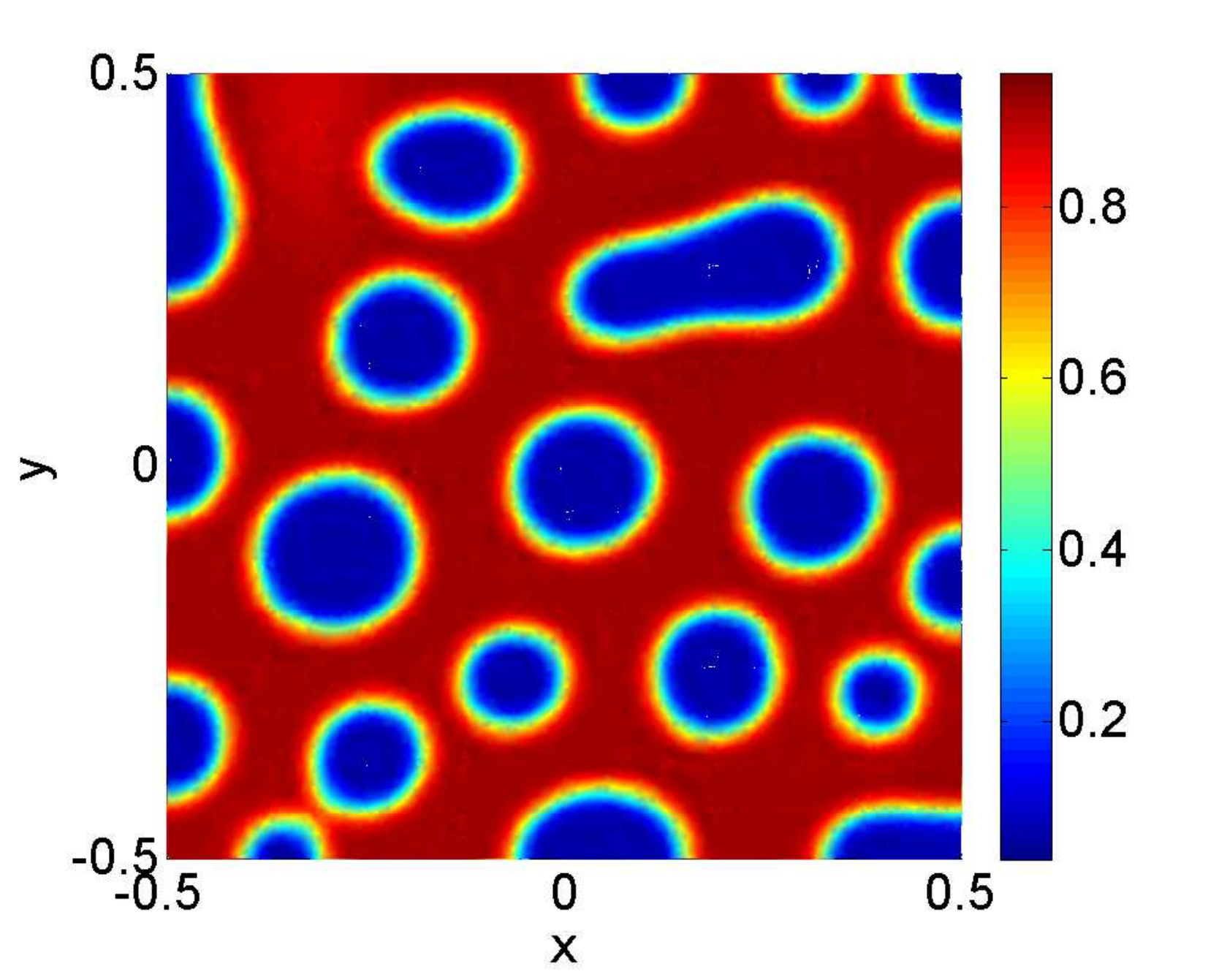}}
\subfigure[$t=0.001$]{\includegraphics[width=0.45\textwidth] {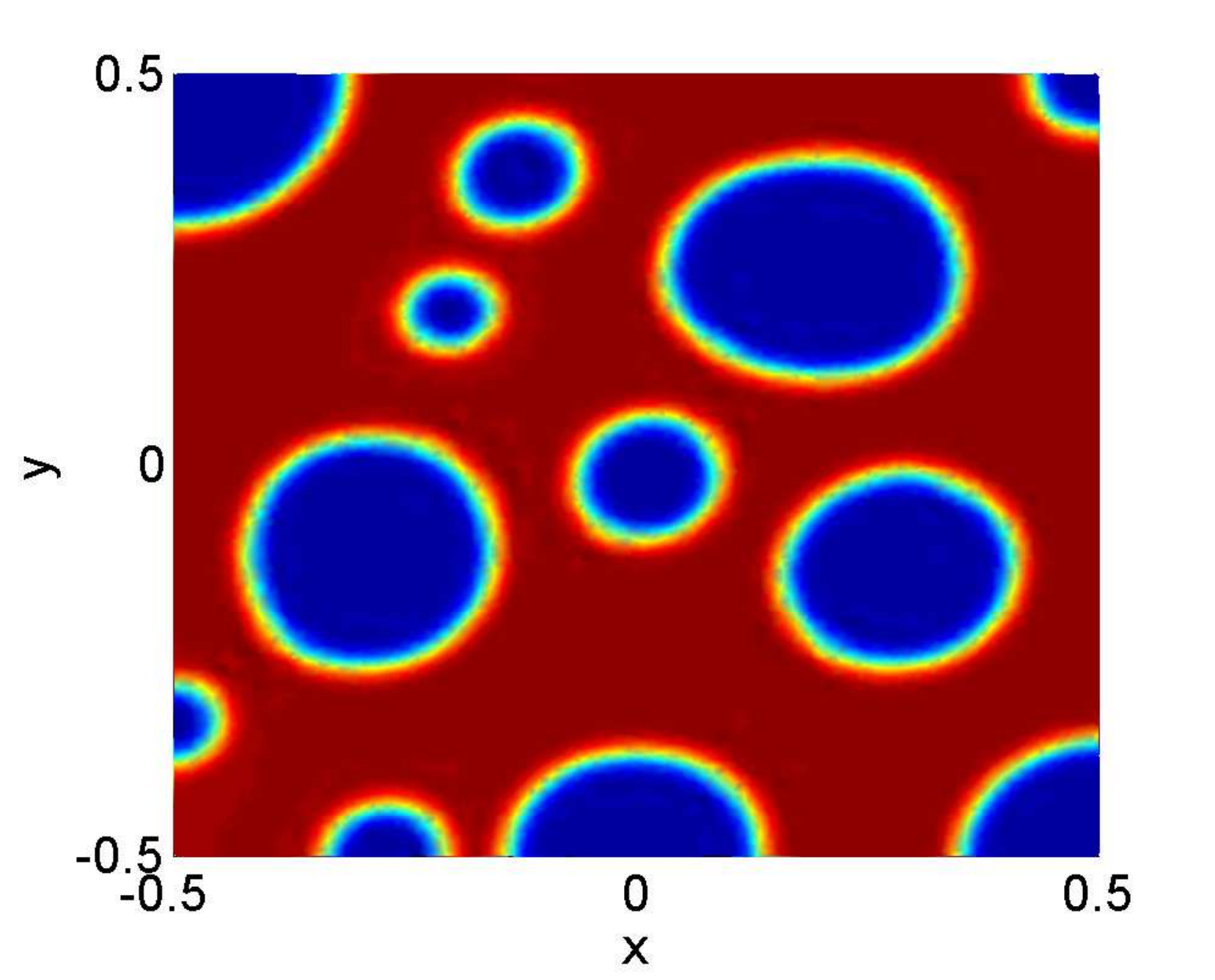}}

\subfigure[]{\includegraphics[width=0.45\textwidth] {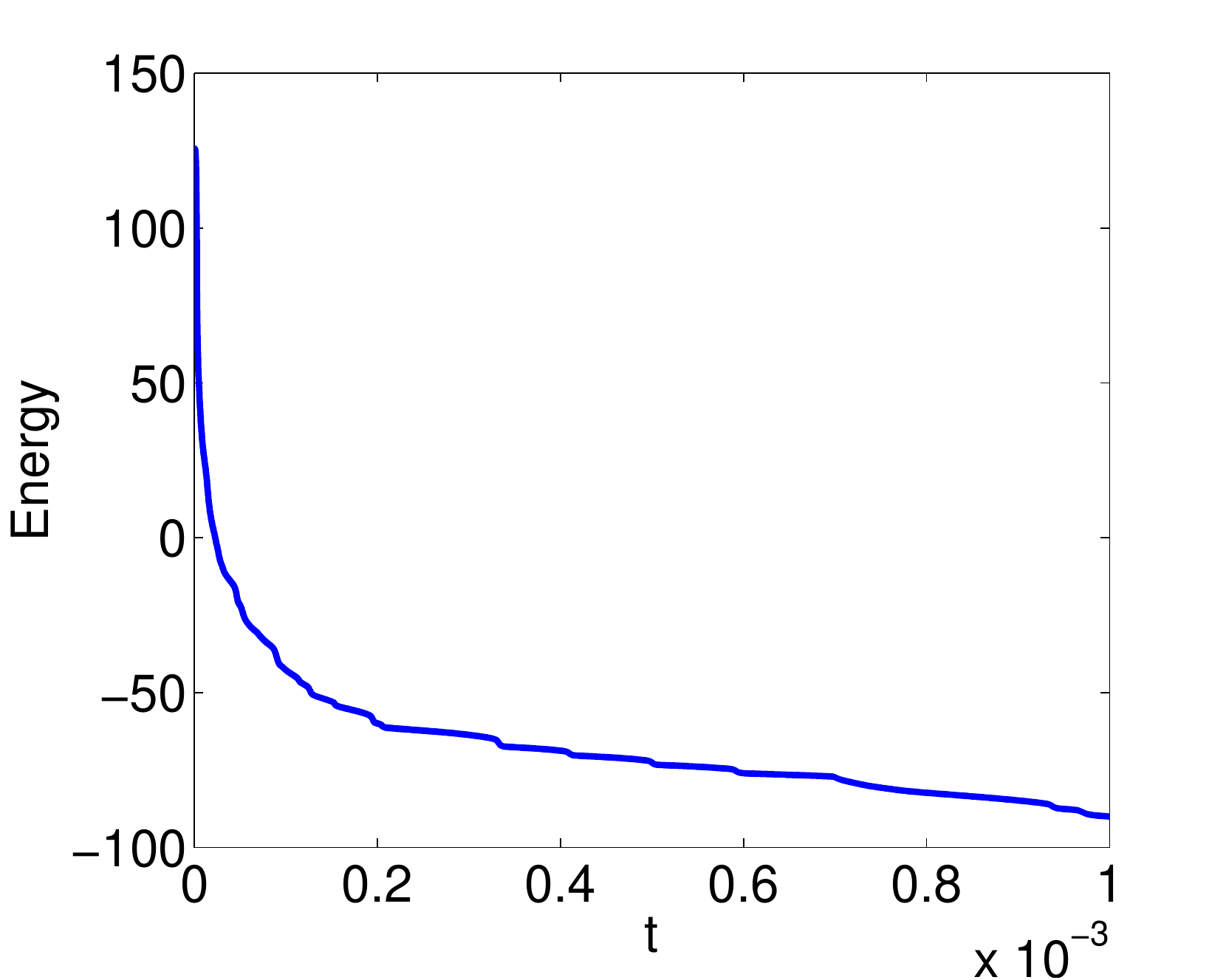}}
\subfigure[]{\includegraphics[width=0.45\textwidth] {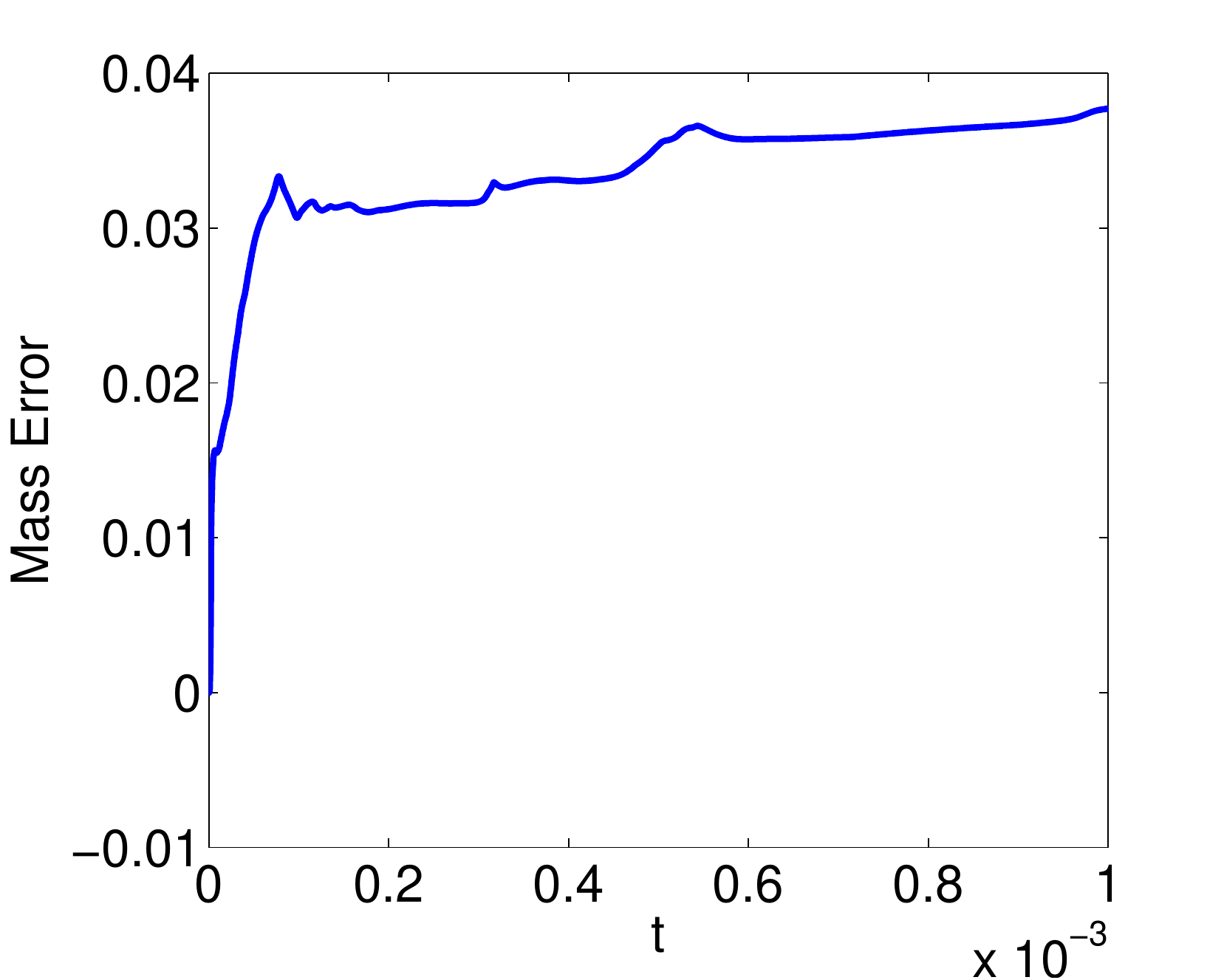}}
\caption{Example \ref{ex4}: (Top-Middle) Solutions with degenerate mobility and  logarithmic potential, (Bottom) energy and mass evolutions. Cubic DG elements and $\Delta t=10^{-7}$.\label{guo_ex7}}
\end{figure}
Fig.~\ref{guo_ex7} shows the evolution of the concentration field. The two phases in the concentration evolution, the phase separation stage and the coarsening process stage can be seen clearly.

\section{Conclusions}\label{conc}
We have presented numerical results for the CH equation under periodic and homogenous Neumann boundary conditions using SIPG discretization in space and AVF method in time. The numerical energy is decreasing in all examples and the numerical results are in good agreement to the those in the literature.

\section*{Acknowledgments}
This work has been supported by Scientific HR Development Program (\"OYP) of the Turkish Higher Education Council (Y\"OK).

\clearpage


\end{document}